\documentclass[12pt]{amsart}
\usepackage{amsmath, amsthm, amssymb}
\usepackage{enumerate}

\theoremstyle{plain}

\newtheorem{thm}{Theorem}[section]
\newtheorem{cor}[thm]{Corollary}
\newtheorem{lem}[thm]{Lemma}
\newtheorem{prop}[thm]{Proposition}

\theoremstyle{definition}
\newtheorem{defi}[thm]{Definition}
\newtheorem{defis}[thm]{Definitions}
\newtheorem{conj}[thm]{Problem}
\newtheorem{conv}[thm]{Convention}
\newtheorem{nota}[thm]{Notation}
\newtheorem{rem}[thm]{Remark}
\newtheorem{rems}[thm]{Remarks}
\newtheorem{exa}[thm]{Example}
\newtheorem{exas}[thm]{Examples}
\newtheorem{sit}[thm]{}

\newcommand{\brem}{\begin{rem}}
\newcommand{\brems}{\begin{rems}}
\newcommand{\erem}{\end{rem}}
\newcommand{\erems}{\end{rems}}
\newcommand{\bexa}{\begin{exa}}
\newcommand{\bexas}{\begin{exas}}
\newcommand{\eexa}{\end{exa}}
\newcommand{\eexas}{\end{exas}}
\newcommand{\bdefi}{\begin{defi}}
\newcommand{\edefi}{\end{defi}}
\newcommand{\bdefis}{\begin{defis}}
\newcommand{\edefis}{\end{defis}}
\newcommand{\bcor}{\begin{cor}}
\newcommand{\ecor}{\end{cor}}
\newcommand{\blem}{\begin{lem}}
\newcommand{\elem}{\end{lem}}
\newcommand{\bconv}{\begin{conv}}
\newcommand{\econv}{\end{conv}}
\newcommand{\bconj}{\begin{conj}}
\newcommand{\econj}{\end{conj}}
\newcommand{\bprop}{\begin{prop}}
\newcommand{\eprop}{\end{prop}}
\newcommand{\bthm}{\begin{thm}}
\newcommand{\ethm}{\end{thm}}
\newcommand{\bnota}{\begin{nota}}
\newcommand{\enota}{\end{nota}}
\newcommand{\bsit}{\begin{sit}}
\newcommand{\esit}{\end{sit}}
\newcommand{\be}{\begin{equation}}
\newcommand{\ee}{\end{equation}}
\newcommand{\bproof}{\begin{proof}}
\newcommand{\eproof}{\end{proof}}
\def\ba{\begin{array}}
\def\ea{\end{array}}
\def\bea{\begin{eqnarray}}
\def\eea{\end{eqnarray}}

\def\bnum{\begin{enumerate}}
\def\enum{\end{enumerate}}

\newcommand{\Spec}{\operatorname{Spec}}

\newcommand{\Frac}{\operatorname{Frac}}

\newcommand{\Proj}{{\operatorname{Proj}}}

\renewcommand{\div}{{\operatorname{div}}}

\newcommand{\la}{\label}

\newcommand{\A}{{\mathbb A}}
\newcommand{\C}{{\mathbb C}}
\newcommand{\Di}{{\mathbb D}}

\newcommand{\N}{{\mathbb N}}
\newcommand{\PP}{{\mathbb P}}

\newcommand{\Q}{{\mathbb Q}}

\newcommand{\V}{{\mathbb V}}
\newcommand{\Z}{{\mathbb Z}}

\newcommand{\G}{{\Gamma}}

\title{Embeddings of $\C^*$-surfaces into weighted projective spaces}

\author{Hubert Flenner}
\address{Fakult\"at f\"ur Mathematik,
Ruhr Universit\"at Bochum, Geb.\ NA 2/72, Universit\"ats\-str.\
150, 44780 Bochum, Germany}
\email{Hubert.Flenner@ruhr-uni-bochum.de}

\author{Shulim Kaliman}
\address{Department of Mathematics,
University of Miami, Coral Gables, FL  33124, U.S.A.}
\email{kaliman@math.miami.edu}

\author{Mikhail Zaidenberg}
\address{Universit\'e
Grenoble I, Institut Fourier, UMR 5582 CNRS-UJF, BP 74, 38402
St.\ Martin d'H\`eres c\'edex, France}
\email{zaidenbe@ujf-grenoble.fr}

\thanks{
\mbox{\hspace{11pt}}{\it 1991 Mathematics Subject
Classification}:
14R05, 14R20.\\
\mbox{\hspace{11pt}}{\it Key words}: 
weighted projective space, $\C^*$-action, $\C_+$-action,
affine surface}

\begin{document}

\begin{abstract}
Let $V$ be a normal affine surface which admits a $\C^*$- and a $\C_+$-action. 
Such surfaces were classified e.g., in \cite{FlZa1, FlZa2}, see also 
the references therein. In this note 
we show that in many cases $V$  can be embedded as a 
principal Zariski open subset into a hypersurface of a 
weighted projective space. In particular, we recover a result of 
D.\ Daigle and P.\ Russell, see Theorem A in \cite{DR}.
\end{abstract}

\maketitle

\section{Introduction}

If $V=\Spec A$ is a normal affine surface 
equipped with an effective $\C^*$-action, then its coordinate ring $A$ 
carries a natural structure of a $\Z$-graded ring
$A=\bigoplus_{i\in \Z} A_i$. 
As was shown in \cite{FlZa1}, such a $\C^*$-action on 
$V$ has a hyperbolic fixed point 
if and only if $C=\Spec A_0$ is a smooth affine curve 
and $A_{\pm 1}\ne 0$.  
In this case the structure of the graded ring $A$ 
can be elegantly described in terms of
a pair $(D_+,D_-)$ of $\Q$-divisors on $C$ with $D_++D_-\le 0$. 
More precisely, $A$ is the graded subring 
$$
A=A_0[D_+,D_-]\subseteq K_0[u,u^{-1}]\,,\quad K_0:=\Frac A_0, 
$$
where for $i\ge 0$
\be\la{pieces}
A_i=\{f\in K_0\mid \div f+iD_+\ge 0\}\,u^i
\quad\mbox{and}\quad
A_{-i}=\{f\in K_0\mid \div f+iD_-\ge 0\}\,u^{-i}\,.
\ee
This presentation  of $A$ (or $V$)
is called in \cite{FlZa1} 
the {\em DPD-presentation}. 
Furthermore two pairs $(D_+,D_-)$ and $(D_+',D_-')$ 
define equivariantly isomorphic surfaces over $C$ 
if and only if they are equivalent that is, 
$$
D_+=D_+'+\div f\quad\mbox{and}\quad D_-=D_-'-\div f 
\quad\mbox{for some }f\in K_0^\times\,.
$$ 
In this note we show that if 
such a surface $V$ admits also a $\C_+$-action then it 
can be $\C^*$-equivariantly embedded (up to normalization) into 
a weighted projective space 
as a hypersurface minus a hyperplane; 
see Theorem \ref{main} 
and Corollary
\ref{k=1} below. In particular 
we recover the following result of 
Daigle and Russell \cite{DR}.

\bthm\la{dr}
Let $V$ be a normal Gizatullin surface\footnote{That is, $V$ admits 
a completion by a linear chain of smooth
rational curves; 
see Section 3 below.} with 
a finite divisor class group. 
Then $V$ can be embedded into a weighted projective plane 
$\PP(a,b,c)$ 
minus a hypersurface. More precisely:
\bnum[(a)]

\item If $V=V_{d,e}$ is toric\footnote{See \ref{last}(a)
below.} 
then $V$ is equivariantly 
isomorphic to the open part\footnote{We use the standard notation 
$\V_+(f)=\{f=0\}$ and $\Di_+(f)=\{f\neq 0\}$.} $\Di_+(z)$ 
of the weighted projective plane $\PP(1,e,d)$ equipped 
with homogeneous coordinates $(x:y:z)$ and with
the $2$-torus action 
$(\lambda_1,\lambda_2).(x:y:z)=(\lambda_1 x:\lambda_2 y: z)$. 

\item If $V$ is non-toric then
$V\cong \Di_+(xy-z^m)\subseteq \PP(a,b, c)$
for some positive integers $a,b,c$ 
satisfying $a+b=cm$ and $\gcd(a,b)=1$.
\enum
\ethm

\section{Embeddings of $\C^*$-surfaces 
into weighted projective spaces}

According to Proposition 4.8 in \cite{FlZa1} every normal affine 
$\C^*$-surface $V$ is equivariantly isomorphic to the normalization of 
a weighted homogeneous surface $V'$ in $\A^4$. In some cases 
(described in {\em loc.cit.}) 
$V'$ 
can be chosen to be a hypersurface in $\A^3$. Cf.\ also \cite{Du} 
for affine embeddings of some other classes of
surfaces.

In Theorem \ref{main} below we show 
that 
any normal $\C^*$-surface $V$ with a $\C_+$-action 
is  the normalization of a principal Zariski open subset of 
some weighted projective hypersurface.

In the proofs we use 
the following observation from \cite{fl}.

\bprop\la{fl}
Let $R=\bigoplus_{i\ge 0}R_i$ be 
a graded $R_0$-algebra of finite type 
containing the field of rational numbers $\Q$. 
If $z\in R_d$, $d>0$, 
is an element of positive degree then 
the group of {\em d}th roots of unity 
$E_d\cong \Z_d$ acts on $R$  
and then also on $R/(z-1)$ via 
$$
\zeta . a =\zeta^i \cdot a\quad\mbox{for}
\quad a\in R_i ,\, \zeta\in E_d,
$$
with ring of invariants $\left(R/(z-1)\right)^{E_d}
\cong \left(R[1/z]\right)_0$. Consequently
$$
(\Spec R/(z-1))/E_d\cong \Di_+(z)
$$
is isomorphic to the complement of the hyperplane 
$\{z=0\}$ in $\Proj(R)$. 
\eprop 

Let us fix the notations.

\bsit\la{nota000}
Let  $V=\Spec A$ be a normal $\C^*$-surface with DPD-presentation
$$A=\C[t][D_+,D_-]\subseteq \C(t)[u,u^{-1}].$$ 
If $V$ carries a $\C_+$-action then according to \cite{FlZa2}, 
after interchanging $(D_+,D_-)$ and passing to an equivalent 
pair, if necessary, we may assume that 
\be\la{dpd}\ba{l}
D_+=-\frac{e_+}{d}[0]\quad\mbox{ with}\quad 0< e_+\le d\,,\\[4pt]
D_-=-\frac{e_-}{d}[0]-\frac{1}{k}D_0
\ea\ee 
with an integral divisor $D_0$, where $D_0(0)=0$.
We choose a polynomial $Q\in \C[t]$ with $D_0=\div( Q)$; so $Q(0)\neq 0$. 
\esit

\bthm\la{main}
Let $F$ be the polynomial
\be\la{pol}
F= x^ky-s^{k(e_++e_-)}Q(s^d/z)z^{\deg Q}
\in \C[x,y,z,s]\,,
\ee
which is weighted homogeneous 
of degree\footnote{We note that $e_++e_-=d(-D_+(0)-D_-(0))\ge 0$.} 
$k(e_++e_-)+d\deg Q$ 
with respect to the weights
\be\la{deg}
\deg x=    e_+ \,,\quad
\deg  y =  ke_-+d\deg Q    \,,\quad
\deg z=    d\,,\quad
\deg s=  1  \,.
\ee
Then the surface $V$ as in \ref{nota000} above 
is equivariantly isomorphic to the normalization of the 
principal Zariski open subset $\Di_+(z)$ of the hypersurface 
$\V_+(F)$
in the weighted projective $3$-space 
\be\la{wps}
\PP=\PP(e_+,ke_-+d\deg Q, d,1)\,.
\ee
\ethm

\bproof
With $s=\sqrt[d]{t}$ the field $L=\Frac(A)[s]$ 
is a cyclic extension of 
$K=\Frac(A)$. Its Galois group is the group of $d$th roots of unity 
$E_d$ acting on $L$ via the identity on $K$ and by 
$\zeta.s=\zeta\cdot s$ 
if $\zeta\in E_d$. Let $A'$ be the normalization of $A$ in $L$. 
According to Proposition 4.12 in \cite{FlZa1}
$$
A'=\C[s][D_+', D_-']\subseteq \C(s)[u,u^{-1}]
$$
with $D'_\pm=\pi_d^*(D_\pm)$, where $\pi_d:\A^1\to\A^1$ is the covering 
$s\mapsto s^d$.
Thus
$$ 
\left(D_+',\, D_-'\right)=\left(-e_+[0],\, -e_-[0]-\frac{1}{k}\pi_d^*(D_0)\right)
=\left(-e_+[0], -e_-[0]-\frac{1}{k}\div (Q(s^d))\right).
$$  
The element $x=s^{e_+}u\in A'_1$ is a generator of $A'_1$ as a $\C[s]$-module. 
According to Example 4.10
in \cite{FlZa1} the graded algebra $A'$ is isomorphic to the normalization of
\be\la{B}
B=\C[x, y,s]/(x^k y-s^{k(e_++e_-)}Q(s^d))\,.
\ee
The cyclic group $E_d$ acts on $A'$ via
$$
\zeta.x=\zeta^{e_+}x\,,\quad 
\zeta.y=\zeta^{ke_-} y\,, \quad
\zeta.s=\zeta s
$$
with invariant ring $A$. Clearly this action stabilizes the subring $B$.  
Assigning to $x, y,z,s$ the degrees as in (\ref{deg}), $F$ as in (\ref{pol}) 
is indeed weighted homogeneous. Since $F(x,y,1,s)=x^k y-s^{k(e_++e_-)}Q(s^d)$, 
the graded algebra 
$$
R=\C]x,y,z,s]/(F)
$$
satisfies $R/(z-1)\cong B$. 
Applying Proposition \ref{fl} $V=\Spec A$ is isomorphic to the normalization of 
$\Di_+(z)\cap \V_+(F)$ in the weighted projective space $\PP$.
\eproof

\brem 
In general not all weights of the weighted projective space $\PP$ in (\ref{wps})
are positive. 
Indeed it can happen that $ke_-+d\deg Q\le 0$. In this case we 
can choose $\alpha\in \N$ with $ke_-+d(\deg Q+\alpha) >0$  and 
consider instead of $F$ the polynomial 
\be\la{pol.111}
\tilde F=x^ky-s^{k(e_++e_-)}Q(s^d/z)z^{\deg Q+\alpha}
\in \C[x,y,z,s]\,,
\ee
which is now weighted homogeneous of degree $k(e_++e_-)+d(\deg Q+\alpha)$ 
with respect to the {\em positive} weights
\be\la{deg.111}
\deg x=    e_+ \,,\quad
\deg  y =  ke_-+d(\deg Q+\alpha)    \,,\quad
\deg z=  d   \,,\quad
\deg s=    1\,.
\ee
As before $V=\Spec A$ is isomorphic to the normalization of the 
principal open subset $\Di_+(z)$ of the hypersurface 
$\V_+(F)$
in the weighted projective space 
$$
\PP=\PP(e_+,ke_-+d(\deg Q+\alpha), d,1)\,.
$$ 
\erem

In certain cases
it is unnecessary  in 
Theorem \ref{main} to pass to normalization. 

\bcor\la{k=1}
Assume that in (\ref{dpd}) one of 
the following conditions is satisfied.
\bnum[(i)]
\item  $k=1$;
\item $e_++e_-=0$, and $D_0$ is a reduced divisor.
\enum
Then $V=\Spec A$ is equivariantly isomorphic to 
the principal open subset $\Di_+(z)$ 
of the weighted projective hypersurface 
$\V_+(F)$ as in (\ref{pol}) 
in the weighted projective space 
$\PP$ from (\ref{wps}). 
\ecor

\bproof
In case (i) the hypersurface in $\A^3$ with equation
$$
F(x,y,1,s)= xy-s^{e_++e_-}Q(s^d)=0
$$
is normal. In other words, 
the quotient $R/(z-1)$ of the graded ring  
$R=\C[x,y,z,s]/(F)$ is normal and so is 
its ring of invariants $\left(R/(z-1)\right)^{E_d}$.
Comparing with Theorem \ref{main} the result follows. 

Similarly, in case (ii) 
$$
F(x,y,1,s)=x^ky-Q(s^d)\,.
$$
Since the divisor $D_0$ is supposed to be reduced and 
$D_0(0)=0$, the  polynomials $Q(t)$ 
and then also $Q(s^d)$ both have  simple roots. 
Hence the  hypersurface $F(x,y,1,s)=0$ in $\A^3$
is again normal, 
and the result follows as before. 
\eproof

\brem\la{smooth} The surface $V$ 
as in \ref{nota000} is smooth if and only if the  divisor $D_0$ is reduced
and $-m_+m_-(D_+(0)+D_-(0))=1$,
where $m_\pm>0$ is the denominator 
in the irreducible representation of $D_\pm(0)$, see
Proposition 4.15 in \cite{FlZa1}. It can happen, however, 
that $V$ is smooth but the surface 
$\V_+(F)\cap \Di_+(z)\subseteq\PP$ has non-isolated singularities. 
For instance, 
if in \ref{nota000}
$D_0=0$ (and so $Q=1$), then $V$ 
is an affine toric surface\footnote{See \ref{last}(a) below.}. 
In fact,
every affine toric surface different from 
$(\A^1_*)^2$ or $\A^1\times\A^1_*$ 
appears in this way, see Lemma 4.2(b) in \cite{FKZ1}. 

In this case the integer $k>0$ can be chosen arbitrarily. 
For any $k>1$, the affine hypersurface 
$V_+(F)\cap \Di_+(z)\subseteq\PP$ 
with equation $x^ky-s^{k(e_++e_-)}=0$ 
has non-isolated singularities
and hence is non-normal. 
Its normalization $V=\Spec A$ can be given 
as the Zariski open part $\Di_+(z)$ of the hypersurface 
$V_+(xy'-s^{e_++e_-})$ in $\PP'=\PP(e_+,e_-,d,1)$ 
(which corresponds to the choice $k=1$). Indeed,
the element $y'=s^{e_++e_-}/x\in K$ with $y'^k=y$ is
integral over $A$. However cf.\ Theorem \ref{dr}(a).
\erem

\bexa\la{D-G} ({\em Danilov-Gizatullin surfaces})
We recall that a Danilov-Gizatullin surface $V(n)$ of index $n$ is the 
complement to a section $S$ in a Hirzebruch surface 
$\Sigma_d$, where $S^2=n>d$. By a remarkable result of 
Danilov and Gizatullin  up to an isomorphism such a surface
 only depends on $n$ 
and neither on $d$ nor on the choice of the section $S$, 
see e.g., \cite{DaGi,  CNR, FKZ3} for a proof.  

According to \cite[\S 5]{FKZ1},  up to conjugation
$V(n)$ 
carries exactly
$(n-1)$ different
$\C^*$-actions. They admit DPD-presentations 
$$
(D_+,D_-)=\left (-\frac{1}{d}[0],\, -\frac{1}{n-d}[1]\right),
\quad\mbox{where}\quad d=1,\ldots,n-1\,.
$$
 Applying 
Theorem \ref{main} with $e_+=1$, $e_-=0$, and $k=n-d$,
the $\C^*$-surface $V(n)$ is the normalization 
of the principal open subset $\Di_+(z)$ 
of the hypersurface $\V_+(F_{n,d})\subseteq \PP(1, d,d,1)$ 
of degree $n$, where
$$
F_{n,d}(x,y,z,s)=x^{n-d}y-s^{n-d}(s^{d}-z)\,.
$$ Taking here $d=1$ it follows that $V(n)$ is isomorphic to 
the normalization of the hypersurface $x^{n-1}y-(s-1)s^{n-1}=0$ 
in $\A^3$. 
\eexa

As our next example, 
let us consider yet another remarkable class 
of surfaces. 
These were studied from different viewpoints in 
\cite[Theorem 1.1]{MM}, 
\cite[Theorem 1.1(iii)]{FlZa3},
\cite[3.8-3.9]{GMMR}, 
\cite[Theorem 1.1. and Example 1]{KK},
\cite[Theorem 1(b) and Lemma 7]{Za}. 
Collecting results from {\em loc.cit.} 
and from this section, 
we obtain the following equivalent characterizations. 

\bthm\la{NT} For a smooth affine surface $V$, 
the following conditions are equivalent. 
\bnum[(i)]

\item $V$ is not Gizatullin and 
admits an effective $\C^*$-action and 
an $\A^1$-fibration $V\to\A^1$ with 
exactly one degenerate fiber,
which is irreducible\footnote{Since $V$ is not Gizatullin there 
is actually a unique 
$\A^1$-fibration  $V\to\A^1$. A surface $V$ as in (i) 
is necessarily a $\Q$-homology plane (or $\Q$-acyclic)
that is, all higher Betti numbers of $V$ vanish.}. 

\item  $V$ is $\Q$-acyclic,  
$\bar k (V)=-\infty$ \footnote{As usual, 
$\bar k$ stands for the logarithmic Kodaira dimension.} 
and $V$
carries a curve 
$\G\cong\A^1$ with $\bar k(V\setminus\G)\ge 0$. 

\item $V$ is $\Q$-acyclic and admits an effective 
$\C^*$- and $\C_+$-actions. Furthermore, 
the $\C^*$-action possesses an orbit closure $\G\cong\A^1$ 
with $\bar k(V\setminus\G)\ge 0$. 

\item The universal cover $\tilde V\to V$ is isomorphic 
to a surface $x^ky-(s^d-1)=0$ in $\A^3$, with the Galois group
$\pi_1(V)\cong E_d$ acting via 
$\zeta.(x,y,s)=(\zeta x,\zeta^{-k}y,\zeta^e s)$,
where $k>1$ and $\gcd (e,d)=1$. 

\item $V$ is isomorphic to the $\C^*$-surface with DPD presentation  
$\Spec \C[t][D_+,D_-]$, where
$$
(D_+,D_-)=\left (-\frac{e}{d}
[0],\,\, \frac{e}{d} [0] -\frac{1}{k} [1]\right)
\qquad\mbox{with}\quad 0<e\le d\quad\mbox{and}\quad k>1\,.
$$

\item $V$ is isomorphic to the Zariski open subset
$$
\Di_+(x^ky-s^d)\subseteq \PP(e,d-ke,1),\qquad\mbox{where}\quad 
0<e\le d\quad\mbox{and}\quad k>1\,. 
$$
\enum
\ethm

\bproof
In view of the references cited above it remains to show 
that the surfaces in $(v)$ and $(vi)$ are isomorphic. 
By Corollary \ref{k=1}(ii) with $e_+=-e_-=e$, the surface
$V$ as in $(v)$ is isomorphic to the principal open subset $\Di_+(z)$ 
in the weighted projective hypersurface
$$
V_+(x^ky-(s^d-z))\subseteq \PP(e,d-ke,d,1)\,.
$$
Eliminating $z$ from the equation $ x^ky-(s^d-z)=0$ yields ($vi$). \eproof

These surfaces admit as well a constructive 
description in terms of 
a blowup process starting from a Hirzebruch
surface, 
see \cite[3.8]{GMMR} and \cite[Example 1]{KK}. 

An affine line $\G\cong\A^1$ on $V$ as in ($ii$) 
is distinguished because it cannot be a fiber 
of any $\A^1$-fibration of $V$. In fact there exists a family 
of such affine lines on $V$, 
see \cite{Za}. 

Some of the surfaces as in Theorem \ref{NT}
can be properly embedded in $\A^3$ as {\em Bertin surfaces} $x^ey-x-s^d=0$, 
see \cite[Example 5.5]{FlZa2} or \cite[Example 1]{Za}.

\section{Gizatullin surfaces with a finite divisor class group}

A 
{\em Gizatullin surface} is a normal affine surface completed by a
zigzag i.e., a linear chain of smooth rational curves. 
By a theorem of Gizatullin \cite{Gi}
such surfaces are characterized by the property that they  
admit two $\C_+$-actions with different general orbits.

In this section we give an alternative proof of the Daigle-Russell
Theorem \ref{dr} 
cited in the Introduction.
It will be deduced from the following result proven in 
\cite[Corollary 5.16]{FKZ2}.

\bprop\la{last} Every normal Gizatullin surface with 
a finite divisor class group is isomorphic to one of the following surfaces. 
\bnum[(a)]
\item The toric surfaces $V_{d,e}=\A^2/E_d$, where   
the group $E_d\cong \Z_d$ of $d$-th roots of unity acts on $\A^2$ via 
$$
\zeta.(x,y)=(\zeta x, \zeta^e y)\,. 
$$

\item The non-toric $\C^*$-surfaces $V=\Spec \C[t][D_+,D_-]$,
where \be\label{dpd00}
\left(D_+,D_-\right)=\left(-\frac{e}{m}[p],\,
\frac{e}{m}[p]-c[q]\right) \quad\mbox{with}\quad c\ge 1,\,\,\,p,q\in
\A^1,\,\,\,p\neq q\,,\ee 
and with coprime integers $e$, $m$ such that $1\le e < m$. 
\enum 

Conversely, any normal affine $\C^*$-surface $V$ as in (a) or (b)  
is a Gizatullin surface with a finite
divisor class group. 
\eprop

Let us now deduce Theorem \ref{dr}. 

\bproof[Proof of Theorem \ref{dr}]
To prove (a), we note that according to \ref{fl} 
the cyclic group $E_d$ acts on the ring 
$\C[x,y,z]/(z-1)\cong\C[x,y]$ 
via $\zeta.x=\zeta x$, $\zeta . y=\zeta^e y$, and $\zeta . z=z$, 
where 
$$
\deg x=1\,,\quad \deg y=e,\quad\mbox{and}\quad \deg z=d\,. 
$$
Hence $\Di_+(z)=\Spec \C[x,y]^{E_d}=V_{d,e}$, 
as required in (a). 

To show (b) we consider
$V=\Spec A$ as in \ref{last}(b), where 
$$
A=\C[t][D_+,D_-]\subseteq \C(t)[u,u^{-1}]\,.  
$$
By definition (\ref{pieces}) the homogeneous pieces
$A_{\pm 1}$ of $A$ are generated as $\C[t]$-modules by the elements
$$
u_+=tu\quad\mbox{and} \quad u_-=(t-1)^cu^{-1}\,,
$$
and similarly $A_{\pm m}$ by  
$$
v_+=t^eu^m\quad\mbox{and} \quad v_-=t^{-e}(t-1)^{cm}u^{-m}\,.
$$ 
Thus
$$
u_+^m=t^{m-e}v_+,\quad u_-^m=t^ev_-,\quad\mbox{and}\quad u_+u_-=t(t-1)^c\,.
$$
The algebra $A$ is the integral closure of the subalgebra 
generated by $u_\pm$, $v_\pm$ and $t$. 

Consider now the normalization $A'$ of $A$ in the field $L=\Frac(A)[u_+']$, 
where 
\be\la{star}
u_+'=\sqrt[d]{v_+}\quad\mbox{with}\quad d=cm\,.
\ee
Clearly the elements  $\sqrt[m]{v_+}=t^{\frac{e-m}{m}}u_+$ 
and then also 
$t^{\frac{e-m}{m}}$ both belong to $L$. Since $e$ and $m$ are coprime 
we can choose $\alpha,\beta\in\Z$ with $\alpha (e-m)+\beta m=1$.
It follows that the element  $\tau:=t^{\frac{1}{m}}=t^{\alpha\frac{e-m}{m}}t^\beta$ 
is as well in $L$ whence being integral over $A$ we have $\tau \in A'$.

The element $u_+'$ as in (\ref{star}) 
also belongs to $A'$ and as well $u_-'=\sqrt[d]{v_-}\in A'$. 
Now $v_+v_-=(t-1)^{cm}$, so taking $d$th roots we get for a 
suitable choice of the root $u_-'$,
\be\la{eq01}
u'_+u'_-=\tau^m -1\,.
\ee
We note that $u_\pm$, $v_\pm$  and $t$ are contained in the subalgebra 
$B=\C[u_+',u_-', \tau]\subseteq A'$. The equation  (\ref{eq01}) 
defines a smooth surface in $\A^3$. Hence $B$ is normal and so 
$$
A'=B\cong\C[u_+',u_-', \tau]/(u_+'u_-'-(\tau^m-1))\,.
$$
By Lemma \ref{04} below, for a suitable $\gamma\in\Z$ 
the integers $a=e-\gamma m$ and $d$ are coprime. 
We may assume as well that $1\le a <d$. 
We let $E_d$ act on $A'$ via $\zeta. u_+'=\zeta^{a} u_+'$ and $\zeta|A={\rm id}_A$. Since $\gcd(a,d)=1$, $A$ is
the invariant ring of this action.   
We claim that the action of $E_d$ on $(u'_+,u_-',\tau)$ is given by 
\be\la{e1}
\zeta. u_+'=\zeta^{a} u_+', \quad \zeta. u_-'=\zeta^{-a} 
u_-'=\zeta^{b}u_-' 
\quad \mbox{and}\quad \zeta.\tau=\zeta^{c}\tau\,, 
\ee 
 where $b=d-a$.
Indeed, the equality
$u_+'^c=t^\frac{e-m}{m}u_+=\tau^{e-m}u_+$  implies
that $\zeta. \tau^{e-m}=\zeta^{ac}\tau^{e-m}$. 
Since $\tau=\tau^{\alpha(e-m)}t^\beta$ the element
$\zeta\in E_d$ acts on $\tau$ via 
$\zeta.\tau=\zeta^{\alpha c a}\tau$. 
In view of the congruence $\alpha a\equiv 1 \mod m$ 
the last expression equals $\zeta^{c}\tau$.
Now the last equality in (\ref{e1}) follows.  
In the equation $u_+'u_-'=\tau^m-1$ the term 
on the right is invariant under $E_d$. 
Hence also the term on the left is. This 
provides the second equality in (\ref{e1}).  

The algebra
$
B=\C[u_+',u_-', \tau]
$
is naturally graded via
$$
\deg u_+'=a, \quad\deg u_-'=b, \quad \mbox{and}\quad \deg \tau= c\,.
$$
According to Proposition \ref{fl} $\Spec A=\Spec A^{\prime E_d}$ 
is the complement of the hypersurface $\V_+(f)$ of degree $d=a+b$ 
in the weighted projective plane
$$
\Proj(B)=\PP(a,b,c),
\quad\mbox{where}\quad  f=u_+'u_-'-\tau^m\,,
$$ 
proving (b).
\eproof

To complete the proof we still have to show the following elementary lemma. 

\blem\la{04}
Assume that $e,m\in\Z$ are coprime. Then for every $c\ge 2$ 
there exists $\gamma\in\Z$ such that $\gamma m-e$ and $c$ are coprime. 
\elem

\bproof
Write $c=c'\gamma$ such that $c'$ and $m$ have no common factor 
and every prime factor of $\gamma$ occurs in $m$. 
Then for every $\gamma\in\Z$ the integers $\gamma m-e$ 
and $\gamma$ have no common 
prime factor. Indeed, such a prime must divide $m$ 
and then also $e=\gamma m-(\gamma m-e)$. 
Hence it is enough to establish the existence of $\gamma\in\Z$ 
such that $\gamma m-e$ and $c'$ are coprime. 
However, the latter is evident since the residue classes 
of $\gamma m$, $\gamma\in \Z$, in $\Z_{c'}$ cover this group.
\eproof

\brem\la{rema}
1. Two triples $(1,e,d)$ and $(1,e',d)$ as in Theorem \ref{dr}(a)
define the same affine toric surface if and only if $ee'\equiv  1\mod d$,
see \cite[Remark 2.5]{FlZa1}. 

2. As follows from Theorem 0.2 in \cite{FKZ2}, the integers $c,m$ 
in Theorem \ref{dr}(b) are invariants of  
the isomorphism type of 
$V$. Indeed, the fractional parts of both divisors $D_\pm$ 
as in (\ref{dpd00}) 
being nonzero
and concentrated at the same point, there is a unique DPD presentation 
for $V$  
up to interchanging $D_+$ and $D_-$, passing to an equivalent
pair and applying an automorphism of the affine line $\A^1=\Spec \C[t]$.

Furthermore, from the proof of Theorem \ref{dr} 
one can easily derive 
that
$$
a\equiv e\mod m
\quad\mbox{and}\quad b=mc-a\equiv -e \mod m\,.
$$
 Therefore  also the pair $(a,b)$ is uniquely determined 
by the isomorphism type of $V$ up to a transposition  
and up to replacing $(a,b)$ by $(a',b')=(a-sm, b+sm)$, while
keeping $\gcd(a',b')=1$. 
\erem

\end{document}